\documentclass[11pt,twoside,leqno]{aomamlt2e} 
  \pageno{971}
\received{June 2, 2004}

 \renewcommand{\thesection}{\arabic{section}}
\makeatletter
\renewcommand{\@seccntformat}[1]{\csname
the#1\endcsname.\hspace{0.5em}\setcounter{Subsec}{0}\setcounter{Subsubsec}{0}}\makeatother

   \newtheorem{theorem}{Theorem}
\newtheorem{corollary}{Corollary}
\newtheorem{lemma}{Lemma}

\newtheorem{thschmidt}{Theorem PV\hskip-4pt}
\newtheorem{thkh}{Khintchine's theorem\hskip-4pt}
\newtheorem{thjar}{Jarn\'{\i}k's theorem\hskip-4pt}

\renewcommand{\Bbb}[1]{\mathbb{#1}}
\newcommand{\N}{{\Bbb N}}         
\newcommand{\I}{{\Bbb I}}
\newcommand{\R}{{\Bbb R}}        
 
\newcommand{\Z}{{\Bbb Z}}         
\newcommand{\K}{{\Bbb K}}         

\newcommand{\cF}{{\cal F}}

\renewcommand{\cH}{{\cal H}}

\newcommand{\cS}{{\cal S}}

\newcommand{\ve}{\varepsilon}

\newcommand{\diam}{r}

\renewcommand{\le}{\leqslant}
\renewcommand{\ge}{\geqslant}

\newcommand{\kgb}{K_{G,B}}
\newcommand{\kgbl}{K_{G',B^{(l)}}}
\newcommand{\kgbf}{K^f_{G,B}}
\newcommand{\hf}{\cH^f}
\newcommand{\hs}{\cH^s}
\newcommand{\p}{\psi}

\begin{document}
\currannalsline{164}{2006} 

 \title{A Mass Transference Principle\\ and  the
Duffin-Schaeffer conjecture\\ for Hausdorff measures}

 \acknowledgements{Research supported by EPSRC  
GR/R90727/01.  

\hglue7pt$^{\textstyle\ast\ast}$Royal Society  
University Research
Fellow}
 
\twoauthors{Victor Beresnevich$^{\textstyle\ast}$}{Sanju Velani\raise2.75pt\hbox{${\textstyle\ast}$}}

 \institution{Institute of Mathematics,
Academy of Sciences of Belarus, Minsk, Belarus
\\\current{University of York, York, England}\\
\email{vb8@york.ac.uk}
\\
\vglue-9pt
University of York, York, England
\\
\email{slv3@york.ac.uk}}

 \shorttitle{A mass transference principle}

 \shortname{Victor Beresnevich and Sanju Velani}

 \centerline{{\it Dedicated to Tatiana Beresnevich}}

\vglue15pt
\centerline{\bf Abstract}
\vglue12pt  A Hausdorff measure version of the
Duffin-Schaeffer conjecture in metric number theory is introduced
and discussed. The general conjecture is  established modulo the
original conjecture.  The key result is a Mass Transference
Principle which allows us to transfer Lebesgue measure theoretic
statements for $\limsup$ subsets of $\R^k$ to Hausdorff measure
theoretic statements. In view of this, the Lebesgue  theory of
$\limsup $ sets  is shown to underpin the general Hausdorff
theory. This is rather surprising since   the latter theory is
viewed to be a subtle refinement of the former.

\section{Introduction}

Throughout $\p:\R^+\to\R^+$ will denote a real, positive
  function 
   and will be referred to as an
\emph{approximating function}. Given an approximating function
$\psi$, a point ${\bf y}=(y_1,\dots,y_k)\in\R^k$ is called {\it
simultaneously $\psi$-approximable} if there are infinitely many
$q\in\N$ and ${\bf p}= (p_1, \ldots,p_k) \in \Z^k$ such that
\begin{equation}
  \left|y_i - \frac{p_i}{q}\right|\ <\
\frac{\psi(q)}{q}  \hspace{9mm}   (p_i,q) = 1 \, , \hspace{5mm} 1
\leq i \leq k \ . \label{1}
\end{equation}
  The set of simultaneously
$\p$-approximable points in $\I^k:=[0,1]^k$ will be denoted by
$\cS_k(\psi)$. For convenience,  we work  within the unit cube
$\I^k$ rather than $\R^k$; it makes full measure results easier
to state and avoids ambiguity. In fact, this is not at all
restrictive as the set of simultaneously $\psi$-approximable
points is invariant under translations by integer vectors.


  The pairwise co-primeness condition
imposed in the above definition clearly ensures that the rational
points $(p_1/q, \ldots,p_k/q)$ are distinct. To some extent the
approximation of points in $\I^k$ by {\em distinct} rational
points should be the main feature when defining $ \cS_k(\psi)$ in
which  case  pairwise co-primeness in (\ref{1}) should be replaced
by the condition that $(p_1, \ldots,p_k,q) =1$. Clearly, both
conditions coincide in the case $k=1$. We shall return to this
discussion  in Section \ref{DSR}.

\Subsec{The Duffin-Schaeffer conjecture\label{tds}}
On making use of the fact that $\cS_k(\psi)$ is a $\limsup$ set, a
simple consequence  of the Borel-Cantelli lemma from probability
theory is that
\begin{equation*}
m(\cS_k(\psi)) = 0  \ \ \ \ \ {\rm if \ } \ \ \ \ \
\sum_{n=1}^{\infty} \left(\phi(n) \, \psi(n)/n \right)^k  \ < \
\infty \ , 
\end{equation*} where $m$ is $k$-dimensional Lebesgue
measure and $\phi$ is the Euler function. In view of this, it is
natural to ask: what happens if the above sum diverges? It is
conjectured that $\cS_k(\psi)$ is of full measure.

\demo{\scshape Conjecture 1}
\begin{equation}
   m(\cS_k(\psi)) = 1  \ \ \ \ \ {\rm if \ } \ \ \ \ \
\sum_{n=1}^{\infty} \left(\phi(n) \, \psi(n)/n \right)^k  \ = \
\infty \ . \label{3}
\end{equation}

When $k=1$, this is the famous Duffin-Schaeffer conjecture in
metric number theory \cite{DS}. Although various partial results
are know,  it remains a major open problem and has attracted much
attention (see \cite{Harman} and references within). For $k \geq
2$, the conjecture was formally stated by Sprind\v{z}uk \cite{Sp}
and settled by Pollington and  Vaughan \cite{PV}.

\begin{thschmidt}
For $k\geq 2${\rm ,} Conjecture {\rm 1} is true.
\end{thschmidt}

If we assume that the approximating function $\psi$ is monotonic,
then we are in good shape thanks to Khintchine's  fundamental
result.

\begin{thkh}
If  $\psi$ is monotonic{\rm ,} then Conjecture $1$ is true.
\end{thkh}

Indeed, the whole point of Conjecture 1  is to remove the
monotonicity condition on $\psi$ from Khintchine's theorem.  Note
that in the case that $\psi$ is monotonic, the
convergence/divergence behavior of the sum in (\ref{3}) is
equivalent to that of $\sum \psi(n)^k$; i.e. the co-primeness
condition imposed in (\ref{1}) is irrelevant.

\Subsec{The Duffin-Schaeffer conjecture for Hausdorff
measures} In this paper, we consider a generalization of
Conjecture 1 which in our view is the `real' problem and the truth
of which yields a complete metric theory. Throughout, $f$ is a
dimension function and $\cH^f$  denotes the Hausdorff $f$-measure;
see Section \ref{HM}. {\em Also, we  assume that $r^{-k} f(r)$ is
monotonic};  this is a natural condition which is  not
particularly restrictive. A straightforward covering argument
making use of the $\limsup$ nature of $\cS_k(\psi)$ implies that
\begin{equation}
\cH^f(\cS_k(\psi)) = 0  \ \ \ \ \ {\rm if \ } \ \ \ \ \
\sum_{n=1}^{\infty}  f( \psi(n)/n ) \; \phi(n)^k  \ < \ \infty \ .
\label{4}
\end{equation}

\noindent In view of this, the following is a `natural'
generalization of Conjecture 1 and can be viewed as the
Duffin-Schaeffer conjecture for Hausdorff measures.

\demo{\scshape Conjecture 2} 
$\displaystyle{\cH^f(\cS_k(\psi)) = \cH^f(\I^k)}$ {\it if}  
$\displaystyle{\sum_{n=1}^{\infty}
f( \psi(n)/n ) \; \phi(n)^k  \ = \ \infty}$.
\Enddemo

Again, in the case that $\psi$ is monotonic  we are in good shape.
This time, thanks to Jarn\'{\i}k's fundamental result.

\begin{thjar}
If  $\psi$ is monotonic{\rm ,} then Conjecture {\rm 2} is true.
\end{thjar}

To be precise, the above  theorem follows on combining
Khintchine's theorem together with  Jarn\'{\i}k's theorem as
stated in \cite[\S8.1]{BDV}; the co-primeness condition imposed
on the set $\cS_k(\psi)$ 
is irrelevant
since $\psi$ is monotonic. The point is that in Jarn\'{\i}k's
original statement, various additional hypotheses on $f$ and
$\psi$  were assumed and they would prevent us from stating the
above clear cut version. Note that Jarn\'{\i}k's theorem together
with (\ref{4}), imply precise Hausdorff dimension results for the
sets $\cS_k(\psi)$;  see \cite[\S1.2]{BDV}.

\Subsec{Statement of results\label{sor}}
Regarding Conjecture 2, nothing seems to be known outside of
Jarn\'{\i}k's theorem which relies on $\psi$ being monotonic. Of
course, the whole point of Conjecture 2 is to remove the
monotonicity condition from Jarn\'{\i}k's theorem. Clearly, on
taking $\cH^f = m $ we have that
$$  Conjecture  \ {\rm  2}  \hspace{4mm}
\Longrightarrow \hspace{4mm}  Conjecture \ {\rm 1}  \ . $$

We shall prove the converse of this statement which turns out to
have obvious but nevertheless  rather unexpected  consequences.

\begin{theorem} 
Conjecture {\rm 1} $ \hspace{4mm} \Longrightarrow \hspace{4mm}  $
Conjecture   {\rm 2}.  \label{main}
\end{theorem}

Theorem \ref{main} together with Theorem PV gives:
\begin{corollary}
For $k\geq 2${\rm ,} Conjecture {\rm 2} is true. \label{cor1}
\end{corollary}

Theorem \ref{main}  gives:
\begin{corollary} Khintchine\/{\rm '}\/s  theorem  $ \hspace{4mm}
\Longrightarrow \hspace{4mm} $  Jarn\'{\i}k\/{\rm '}\/s  theorem.
\label{cor2}
\end{corollary}

It is remarkable that Conjecture 1,  which is only concerned with
the metric theory of $\cS_k(\psi)$ with respect to the ambient
measure $m$,  underpins the whole  general  metric theory. In
particular, as a consequence of Corollary \ref{cor2}, if $\psi$ is
monotonic then Hausdorff dimension results for $\cS_k(\psi)$ (i.e.\
the general form of the Jarn\'{\i}k-Besicovitch theorem) can in
fact be obtained via Khintchine's Theorem. At first, this seems
rather counterintuitive. In fact, the dimension results for
monotonic $\psi$ are a trivial consequence of Dirichlet's theorem
(see \S\ref{secJB}).

The key to establishing Theorem \ref{main} is the Mass
Transference Principle of Section \ref{secMTP}. In short, this allows us
to transfer $m$-measure theoretic statements for $\limsup$ subsets
of $\R^k$ to $\cH^f$-measure theoretic statements. In
Section \ref{secmtpg}, we state a general  Mass Transference Principle
which allows us to obtain the analogue of Theorem \ref{main} for
$\limsup$ subsets of  locally compact metric  spaces.

\section{Preliminaries}

Throughout $(X,d)$ is a metric space such that for every $\rho>0$
the space $X$ can be covered by  a  countable collection of balls
with diameters $<\rho$. A ball $B=B(x,r):=\{y\in X:d(x,y)\le r\}$
is defined by a fixed centre and radius, although these in general
are not uniquely determined by $B$ as a set. By definition, $B$ is
a subset of $X$. For any $\lambda>0$, we denote by  $\lambda
  B$  the ball $B$ scaled by a factor $\lambda$; i.e.  $\lambda B(x,
  r):= B(x, \lambda r)$.

\Subsec{Hausdorff measures\label{HM}}
In this section we give a brief account of Hausdorff measures. A
{\em dimension function} $f \, : \, \R^+ \to \R^+ $ is a continuous,
nondecreasing function such that $f(r)\to 0$ as $r\to 0 \, $. Given
a ball $B=B(x,r)$, the quantity
\begin{equation}\label{e:004}
V^f(B)\,:=\,f(r)
\end{equation}
will be referred to as the  {\em $f$-volume of $B$}.  If $B$ is a
ball in $\R^k$, $m$ is $k$-dimensional Lebesgue measure and
$f(x)=m(B(0,1))x^k$, then $V^f$ is simply the volume of $B$ in the
usual geometric sense; i.e. $V^f(B)=m(B)$. In the case when
$f(x)=x^s$ for some $s\geq0$, we write $V^s$ for $V^f$.

The Hausdorff $f$-measure with respect to the dimension function
$f$ will be denoted throughout by ${\cal H}^{f}$ and is defined as
follows. Suppose $F$ is  a  subset  of $(X,d)$. For $\rho
> 0$, a countable collection $ \left\{B_{i} \right\} $ of balls in
$X$ with $\diam(B_i) \leq \rho $ for each $i$ such that $F \subset
\bigcup_{i} B_{i} $ is called a {\em $ \rho $-cover for $F$}.
Clearly such a cover exists for every $\rho > 0$. For a dimension
function $f$ define $$
  {\cal H}^{f}_{\rho} (F) \, = \, \inf \ \sum_{i} V^f(B_i),
$$
where the infimum is taken over all $\rho$-covers of $F$. The {\it
Hausdorff $f$-measure} $ {\cal H}^{f} (F)$ of $F$ with respect to
the dimension function $f$ is defined by   $$ {\cal H}^{f} (F) :=
\lim_{ \rho \rightarrow 0} {\cal H}^{f}_{\rho} (F) \; = \;
\sup_{\rho > 0 } {\cal H}^{f}_{\rho} (F) \; . $$

A simple consequence of the definition of $ {\cal H}^f $ is the
following useful fact.

\begin{lemma}
If $ \, f$ and $g$ are two dimension functions such that the ratio
$f(r)/g(r) \to 0 $ as $ r \to 0 ${\rm ,} then ${\cal H}^{f} (F) =0 $
whenever ${\cal H}^{g} (F) < \infty $. \label{dimfunlemma}
\end{lemma}

In the case that  $f(r) = r^s$ ($s \geq 0$), the measure $ \hf $
is the usual $s$-dimensional Hausdorff measure $\hs $ and the
Hausdorff dimension $\dim F$ of a set $F$ is defined by $$ \dim \,
F \, := \, \inf \left\{ s : {\cal H}^{s} (F) =0 \right\} = \sup
\left\{ s : {\cal H}^{s} (F) = \infty \right\} . $$ In particular
when $s$ is an integer and $X = \R^s$, $\hs$ is comparable to
the $s$-dimensional Lebesgue measure. Actually, $\hs$ is a constant
multiple of the $s$-dimensional Lebesgue measure but we shall not
need this stronger statement.



For further details see \cite{falc, mat}. A general and classical
method for obtaining a lower bound for the Hausdorff $f$-measure
of an arbitrary set $F$ is the following mass distribution
principle.

\demo{\scshape Lemma (Mass Distribution Principle)} {\em
  Let $ \mu $ be a probability measure supported on a subset $F$ of $  
(X,d) $.
Suppose there are  positive constants $c$ and $r_o$ such that} $$
\mu ( B ) \leq \, c \;  V^f(B) \;
$$ {\em for any ball $B$ with
radius $r \leq r_o \, $.  If $E$ is a subset of $F$ with $\mu(E) =
\lambda > 0$  then $ {\cal H}^{f} (E) \geq \lambda/c \, $. }

\Proof     If $ \left\{B_{i} \right\} $ is a
$\rho$-cover of $E$ with $\rho \leq r_o$ then $$ \lambda = \mu(E)
= \mu \left( \cup_{i} B_i \right) \leq \sum_i \mu \left(  B_i
\right) \leq c \sum_i     V^f(B_i)  \; . $$ It follows that $
{\cal H}^{f}_{\rho } (E) \geq \lambda/c $ for any $\rho \leq r_o$.
On letting $\rho \to 0 \, $,  the quantity ${\cal H}^{f}_{ \rho }
(E)$ increases and so we obtain the required result.
\Endproof\vskip4pt


The following basic covering lemma will be required at various
stages \cite{jh},~\cite{mat}.

\begin{lemma}[The $5r$ covering lemma]\label{5r}
Every family ${\cal F}$ of balls of uniformly bounded diameter in
a metric space $(X,d)$ contains a disjoint subfamily ${\cal G} $
such that $$ \bigcup_{B \in {\cal F} } B \ \subset \ \bigcup_{B
\in {\cal G} } 5B \ \ \ . $$
\end{lemma}


\Subsec{Positive and full measure sets}
  Let $\mu$ be a finite measure supported on
$(X,d)$. The measure $\mu$ is said to be {\em doubling} if there
exists a  constant $\lambda
> 1 $ such that for $x \in X$ $$\mu(B(x,2r)) \, \leq \, \lambda\,   
  \mu(B(x,r)) \ . $$
Clearly, the measure $\cH^k$ is a doubling measure on $\R^k$.  In
this section we state two measure theoretic results which will be
required during the course of the paper.

\begin{lemma}\label{lem1a}
Let $(X,d)$ be a metric space and let $\mu$ be a finite doubling
measure on $X$ such that any open set is $\mu$ measurable. Let $E$
be a Borel subset of\/ $X$. Assume that there are constants
$r_0,c>0$ such that for any ball $B$ with $r(B)<r_0$  and center
in $X${\rm ,}  we have that $\mu(E\cap B)\ge c\ \mu(B)$. Then{\rm ,} for any
ball $B$
$$\mu(E\cap B) \, =  \, \mu(B)  \ . $$
\end{lemma}

\begin{lemma}\label{lem3a}
Let $(X,d)$ be a metric space and $\mu$ be a finite measure on $X$.
Let $B$ be a ball in $X$ and $E_n $ a sequence of
$\mu$-measurable sets. Suppose there exists a constant $c > 0$
such that $\limsup_{n \to \infty} \mu(B \, \cap \, E_n) \ge c \;
\mu(B)$. Then $$ \mu (B \, \cap \, \limsup_{n\to\infty}E_n ) \ \ge
\ c^2 \, \mu(B) \ \ . $$
\end{lemma}

For the details regarding these two lemmas see \cite[\S8]{BDV}.

\section{A mass transference principle \label{secMTP}}

Given a  dimension function $f$ and a ball $B=B(x,r)$ in $\R^k$,
we define another ball
\begin{equation}\label{e:006}
\textstyle B^f:=B(x,f(r)^{1/k}) \ .
\end{equation}
When $f(x)=x^s$ for some $s>0$ we also adopt the notation $B^s$,
i.e.\ $ B^s:=B^{(x\mapsto x^s)}. $ It is readily verified that
\begin{equation}\label{e:008}
B^k=B.
\end{equation}

Next, given a collection $K$ of balls in $\R^k$, denote by $K^f$
  the collection of balls obtained from $K$ under the
transformation (\ref{e:006}); i.e. $K^f := \{ B^f : B \in K \} $.

The following property immediately follows from (\ref{e:004}),
(\ref{e:006}) and (\ref{e:008}):
\begin{equation}\label{e:009}
V^k(B^f)=V^f(B^k)\qquad\text{for any ball $B$.}
\end{equation}
Note that (\ref{e:009}) could have been taken to be a definition
in which case (\ref{e:006}) would follow.

Recall that $\cH^k$ is comparable to the $k$-dimensional Lebesgue
measure $m$. Trivially, for any ball $B$ we have that $V^k(B)$ is
comparable to $m(B)$. Thus there are constants
$0<c_1<1<c_2<\infty$ such that for any ball $B$
\begin{equation}\label{e:010}
c_1\ V^k(B)\le \cH^k(B)\le c_2\ V^k(B).
\end{equation}

\noindent In fact, we have the stronger statement that $ \cH^k(B)
$ is a constant multiple of $V^k(B) $. However, the analogue of
this stronger statement is not necessarily true in the general
framework considered in Section \ref{secmtpg} whereas (\ref{e:010}) is.
Therefore, we have opted to work with (\ref{e:010}) even in our
current setup. Given a sequence of balls $B_i$, $i=1,2,3,\ldots$,
as usual its limsup set is
$$
\limsup_{i\to\infty}B_i:=\bigcap_{j=1}^\infty\ \bigcup_{i\ge j}B_i
\  .
$$

\noindent The following theorem is without doubt  the main result
of this paper.  It is the  key to establishing the
Duffin-Schaeffer conjecture for Hausdorff measures.

\begin{theorem}[Mass Transference Principle]\label{thm3}
Let $\{B_i\}_{i\in\N}$ be a sequence of balls in $\R^k$ with
$\diam(B_i)\to 0$ as $i\to\infty$. Let $f$ be a dimension function
such that $x^{-k}f(x)$ is monotonic and  suppose that for  any ball
$B$ in $\R^k$
\begin{equation}\label{e:011}
\cH^k\big(\/B\cap\limsup_{i\to\infty}B^f_i{}\,\big)=\cH^k(B) \ .
\end{equation}
Then{\rm ,} for any ball $B$ in $\R^k$
\begin{equation*}
\cH^f\big(\/B\cap\limsup_{i\to\infty}B^k_i\,\big)=\cH^f(B) \ .
\end{equation*}
\end{theorem}

 {\em Remark}\/ 1. $\cH^k$ is comparable to the Lebesgue
measure $m$  in $\R^k$. Thus (\ref{e:011}) simply states that the
set $\limsup B_i^f$ is of  full $m$ measure  in $\R^k$, i.e.\ its
complement in $\R^k$ is of $m$ measure zero.

 \demo{Remark    $2$}  In the statement of Theorem~\ref{thm3}
the condition $\diam(B_i)\to0$ as $i\to\infty$ is redundant.
However,
  it is included to avoid unnecessary further discussion.

\demo{Remark  $3$} If $x^{-k}f(x) \to l $ as $x \to 0 $
and $l$ is finite then the above statement is relatively
straightforward to establish. The main substance of the Mass
Transference Principle is when $x^{-k}f(x)\to\infty$ as $x\to0$.
In this case, it trivially follows via Lemma \ref{dimfunlemma}
that $\cH^f(B)=\infty$.

\Subsec{Proof of Theorem {\rm \ref{main}}\label{thmmain}}
First of all let us dispose of the case that $\psi(r)/r
\nrightarrow 0 $ as $r \to \infty$. Then trivially, $\cS_k(\psi)=
\I^k$ and the result is obvious. Without loss of generality,
assume that $\psi(r)/r \to 0 $ as $r \to \infty$. We are given
that $ \sum f( \psi(n)/n ) \; \phi(n)^k =  \infty $. Let $
\theta(r) := r \, f( \psi(r)/r )^{1/k}$. Then $\theta$ is an
approximating function and $ \sum (\phi(n) \, \theta(n)/n )^k  =
\infty $. Thus, on using the supremum norm, Conjecture 1 implies
that $ \cH^k(B \cap \cS_k(\theta)) = \cH^k(B \cap \I^k)$ for any
ball $B$ in $\R^k$. It now follows via the Mass Transference
Principle that $ \cH^f(\cS_k(\psi)) = \cH^f(\I^k) $ and this
completes the proof of Theorem \ref{main}.

\Subsec{The Jarn\'{\i}k-Besicovitch theorem}\label{secJB}
In the case $k=1$ and $\psi( x) := x^{-\tau}\!$, let us write
$\cS(\tau)$ for $\cS_k(\psi)$.  The Jarn\'{\i}k-Besicovitch
theorem states that $\dim \cS(\tau) = d:= 2/(1+\tau) $ for $\tau
> 1$. This fundamental result is easily deduced on combining
Dirichlet's theorem with the Mass Transference Principle.

Dirichlet's theorem states that for any irrational $y \in \R$,
there exists infintely many reduced rationals $p/q$ ($q>0$) such
that $|y - p/q| \leq q^{-2}$. With $f(x) := x^{d}$, (\ref{e:011})
is trivially satisfied  and the Mass Transference Principle
implies that $\cH^{d} (\cS(\tau) ) = \infty $. Hence $\dim
\cS(\tau) \geq d $. The upper bound is trivial. Note
that we have actually proved a lot more than simply the
Jarn\'{\i}k-Besicovitch theorem. We have proved that the
$s$-dimensional  Hausdorff measure $\cH^{s}$ of $\cS(\tau)$ at
the critical exponent $s=d$ is infinite.

\section{The $K_{G,B} $  covering lemma}

Before establishing the Mass Transference Principle we state and
prove the following covering lemma, which provides an equivalent
description of the full measure property (\ref{e:011}).

\begin{lemma}[The $K_{G,B} $ lemma] \label{lem1}
Let $\{B_i\}_{i\in\N}$ be a sequence of balls in $\R^k$ with
$\diam(B_i)\to 0$ as $i\to\infty$. Let $f$ be a dimension function
  and for any  ball $B$
in $\R^k$ suppose that {\rm (\ref{e:011})} is satisfied.  Then for any
  $B$  and any $G>1$ there is a finite
sub-collection  $\kgb\subset\{B_i\,:\,i\ge G\}$ such that the
corresponding  balls in $\kgbf$ are disjoint{\rm ,} lie inside $B$ and
\begin{equation}\label{e:014}
\cH^k\Big(\bigcup^\circ_{L\in \kgbf}L\Big) \ \ge \ \kappa
  \ \cH^k(B) \,
  \hspace{12mm}\text{with}\quad\kappa := \mbox{\footnotesize$  
\frac{1}{2} $}
  (\mbox{\footnotesize$ \frac{c_1}{c_2} $} )^2 10^{-k} \ \ .
\end{equation}
\end{lemma}

 {\it Proof of Lemma}~\ref{lem1}. Let $\cF:=\{B^f_i\,:\,
B^f_i\ \cap \mbox{\footnotesize{$\frac{1}{2}$}}B \neq \emptyset \,
,\ i\ge G\}$. Since, $f(x)\to0$ as $x\to0$ and $\diam(B_i)\to0$ as
$i\to\infty$ we can ensure that every ball in $ \cF$ is contained
in $ B $ for $i$ sufficiently large. In view of the $5r$ covering
lemma (Lemma \ref{5r}), there exists a disjoint sub-family ${\cal
G} $ such that $$ \bigcup_{B^f_i\in {\cal F} } B^f_i \ \subset \
\bigcup_{B^f_i \in {\cal G} } 5B^f_i \ \ \ . $$ It follows that
\begin{eqnarray*}
{\cal H}^k \left(\bigcup_{B^f_i \in {\cal G} } 5B^f_i\right)  \
\geq \ {\cal H}^k \big( \mbox{\footnotesize{$\frac{1}{2}$}} B \cap
\limsup_{i\to\infty}B^f_i \big)  \stackrel{(\ref{e:011})}{\ = \ }
   \cH^k(\mbox{\footnotesize{$\frac{1}{2}$}} B \big)  
\stackrel{(\ref{e:010})}{\ \geq \ }
\frac{c_1}{c_2}  \ 2^{-k} \; \cH^k(B) \ .
\end{eqnarray*}

\noindent However, since ${\cal G} $ is a disjoint collection of
balls we have that
\begin{eqnarray*}
{\cal H}^k \left(\bigcup_{B^f_i \in {\cal G} } 5B^f_i\right)
\stackrel{(\ref{e:010})}{\ \leq \ } \frac{c_2}{c_1}  \ 5^{k} \ \
{\cal H}^k \left(\bigcup^\circ_{B^f_i \in {\cal G} } B^f_i\right)
\ .
\end{eqnarray*}
Thus, \begin{equation}
  {\cal H}^k \left(\bigcup^\circ_{B^f_i \in
{\cal G} } B^f_i\right) \ \geq \ \left(\frac{c_1}{c_2}\right)^2 \
10^{-k} \; \ \cH^k(B) \ . \label{kgb001} \end{equation}  The balls
$B^f_i \in {\cal G}$ are disjoint, and since $\diam(B^f_i)\to0$ as
$i\to\infty$ we have that $$ {\cal H}^k \left(\bigcup^\circ_{B^f_i
\in {\cal G}\, : \, i \geq j  } B^f_i\right) \ \to \  0
\hspace{8mm} {\rm as } \hspace{4mm} j \to \infty \  \ . $$ Thus,
there exists some $j_0 > G $ for which
\begin{equation}  {\cal H}^k \left(\bigcup^\circ_{B^f_i \in {\cal G} \,  
: \, i
\geq j_0  } B^f_i\right) \ < \ \frac12 \;
\left(\frac{c_1}{c_2}\right)^2 \ 10^{-k} \; \ \cH^k(B) \ .
\label{kgb002} \end{equation}  Now let $\kgb := \{B_i : B^f_i \in
{\cal G} , i < j_0 \} $. Clearly, this is a finite sub-collection
of $\{B_i\,:\,i\ge G\}$. Moreover, in view of (\ref{kgb001}) and
(\ref{kgb002}) the collection $\kgbf$
satisfies the desired properties.
\Endproof\vskip4pt

Lemma~\ref{lem1} shows that the full measure property
(\ref{e:011}) of the Mass Transference Principle implies the
existence of the collection $\kgbf$ satisfying  (\ref{e:014}) of
the $\kgb$ Lemma. For completeness, we prove  that the converse is also
true.

\begin{lemma}\label{lem1-}
Let $\{B_i\}_{i\in\N}$ be a sequence of balls in $\R^k$ with
$\diam(B_i)\to 0$ as $i\to\infty$. Let $f$ be a dimension function
and for any   ball $B$  and any $G>1${\rm ,} assume that there is a
collection  $\kgbf$  of balls satisfying {\rm (\ref{e:014})} of
Lemma~{\rm \ref{lem1}.}
  Then{\rm ,} for any ball $B$ the full measure property {\rm (\ref{e:011})} of the
  Mass Transference Principle is
satisfied.
\end{lemma}

{\it Proof of Lemma~{\rm \ref{lem1-}}}. For any ball $B$ and
any $G\in\N$, the collection  $\kgbf$ is contained in $B$ and is a
finite sub-collection of $\{B^f_i\}$ with $i\ge G$.  We define
$$E_G \, := \, \bigcup_{L\in\kgbf}L \ .$$ Since $\kgbf$ is finite, we
have that \begin{equation*} \limsup_{G\to\infty}E_G \ \subset \
B\cap\limsup_{i\to\infty}B^f_i \ . \label{-inc} \end{equation*} It
follows from (\ref{e:014}) that $ \cH^k(E_G)  \ge  \kappa \,
\cH^k(B) $ which together with Lemma~\ref{lem3a} implies that
  $\cH^k(\limsup_{G\to\infty} E_G)\ge
\kappa^2\,\cH^k(B)$.  
Hence, $\cH^k(B\cap\limsup_{i\to\infty} B^f_i)\ge
\kappa^2\,\cH^k(B)$. The measure $\cH^k$ is  doubling and so the
statement of the lemma  follows on applying Lemma~\ref{lem1a}.
\Endproof 

In short, Lemmas \ref{lem1} and \ref{lem1-} establish the
equivalence: $  \mbox{ (\ref{e:011})} \iff \mbox{(\ref{e:014})} $.


\vglue-22pt
\phantom{up}
\section{Proof of Theorem \ref{thm3} (Mass Transference Principle)}
\vglue-8pt

We start by  considering the case that $x^{-k}f(x) \to l $ as $x
\to 0 $ and $l$ is finite. If $l=0$, then Lemma \ref{dimfunlemma}
implies that $\cH^f(B)=0$ and since $ B\cap\limsup B^k_i \subset
B$ the result follows. If $l \neq 0 $ and is finite then  $\cH^f$
is comparable to $\cH^k$ (in fact, $\cH^f = l \,  \cH^k$).
Therefore the required statement follows on showing that
$\cH^k\big(\/B\cap\limsup_{i\to\infty}B^k_i\,\big)=\cH^k(B)$. This
can be established by first noting that the ratio of the radii of
the balls $B^k_i$ and $B^f_i$ are uniformly bounded between
positive constants and then
  adapting the proof of Lemma~\ref{lem1-} in the obvious manner.

In view of the above discussion, we can assume without loss of
generality that $$x^{-k}f(x) \ \to \ \infty \hspace{6mm} {\rm as }
\hspace{6mm} x\to0 \ \ . $$ Note that in this case, it trivially
follows via Lemma \ref{dimfunlemma} that $\cH^f(B)=\infty$. Fix
some arbitrary bounded ball $B_0$ of $\R^k$. The statement of the
Mass Transference Principle will therefore follow on showing that
$$\cH^f (B_0 \cap \limsup B_i) = \infty  \ . $$

\noindent To achieve this we proceed as follows. For any constant
$\eta>1$, our aim is to construct a Cantor subset $\K_\eta$ of
$B_0 \cap \limsup B_i$ and a probability measure $\mu$ supported
on $\K_{\eta}$ satisfying the condition that for an arbitrary ball
$A$ of sufficiently small radius $r(A)
$
\begin{equation}
\mu(A) \; \ll \; \frac{V^f(A)}{\eta}  \;   , \label{task}
\end{equation}
where the implied constant in the Vinogradov symbol ($\ll$) is
absolute. By the Mass Distribution Principle, the above inequality
implies that $$ {\cal H}^f(\K_{\eta}) \; \gg \; \eta \;. $$ Since
$ \K_{\eta}\subset B_0 \cap \limsup B_i$, we obtain that ${\cal
H}^f \left( B_0 \cap \limsup B_i \right)\gg \eta$. However, $\eta
$ can be made arbitrarily large whence ${\cal H}^f \left( B_0 \cap
\limsup B_i \right)=\infty$ and this proves Theorem \ref{thm3}.

In view of the above outline, the whole strategy of our proof is
centred around the  construction of a `right type' of Cantor set
$\K_{\eta}$ which  supports a measure $\mu$ with the desired
property.

\Subsec{The desired properties  of $\K_{\eta}$}
In this section we summarize the  desired  properties of  the
Cantor set $\K_\eta$. The existence of  $\K_\eta$  will be
established in the next section. Let $$
\K_\eta:=\bigcap_{n=1}^\infty \K(n)  \ , $$ where each {\em level}
$\K(n)$ is a finite union of  disjoint balls such that
\begin{equation*}\label{e:016}
   \K(1)\supset \K(2) \supset \K(3) \supset \ldots \ \ .
\end{equation*}
Thus, the levels are nested. Moreover, if $K(n)$  denotes the
collection of balls which  constitute level $n$,  then $K(n)
\subset \{B_i : i \in \N \} $ for each $n \geq 2$.  We will define
$K(1) := B_0$.  It is then clear that $\K_\eta$ is a subset  of
$B_0 \cap \limsup B_i$. It will be convenient to also refer to the
collection $K(n)$ as  the $n$-th level. Strictly speaking, $\K(n)=
\bigcup_{ B \in K(n)} B  $ is the $n$-th level. However, from the
context it will be clear what we mean and no ambiguity should
arise.

The  construction is inductive and the  general idea is as
follows. Suppose the $(n-1)$-th   level $\K(n-1)$  has been
constructed. The next level is constructed by `looking' locally at
each ball from the previous level.  More precisely,  for every
ball $B\in K(n-1)$ we construct the $(n,B)$-{\em local level}
denoted by $K(n,B)$ consisting  of balls contained in $B$. Thus
\begin{equation*}\label{e:017}
K(n)\ := \ \bigcup_{B\in K(n-1)}K(n,B)\qquad\text{and}\qquad
\K(n)\ := \ \bigcup_{B\in K(n-1)}\K(n,B) \ \ ,
\end{equation*}
where $$\K(n,B)\ := \ \bigcup_{L\in K(n,B)} \!\!\!\! L \ = \ B\cap
\K(n) \ \ . $$

\noindent As  mentioned above, the balls in each level will be
disjoint. Moreover, we ensure that  balls in each level  scaled by
a factor of three are  disjoint. This is property (P1) below. This
alone is not  sufficient to obtain the required  lower bound for
$\cH^f(\K_\eta)$.   For this purpose, every local level will be
defined as a union of {\em local sub-levels}. The $(n,B)$-local
level will take on  the following form
\begin{equation*}\label{e:018}
K(n,B):=\bigcup_{i=1}^{l_B}K(n,B,i) \ ,
\end{equation*}
where $l_B$ is the number of local sub-levels (see property (P5)
below) and $K(n,B,i)$ is the $i$-th local sub-level. Within each
local sub-level $K(n,B,i)$, the separation of balls is much more
demanding than simply property (P1) and is given by property (P2)
below.

To achieve our main objective,   the lower bound for
$\cH^f(\K_\eta)$, we will  require a controlled build up of
`mass' on the balls in every sub-level. The mass is related to the
$f$-volume $V^f$  of the balls in the construction and the overall
number of sub-levels. These are governed by properties (P3) and
(P5) below.

Finally, we will require that the $f$-volume of  balls from one
sub-level to the next decreases sufficiently fast. This is
property (P4) below. However,  the total $f$-volume within any one
sub-level remains about the same. This is a consequence of
property (P3) below.

We now formally state the properties (P1)--(P5)  discussed above
together with a trivial property (P0).

\smallbreak {\it The properties of levels and sub-levels of $\K_\eta$}
\begin{enumerate}

\item[{\bf(P0)}] $K(1)$ consists of one ball, namely $B_0$.

\item[{\bf(P1)}] For any $n\ge 2$ and any $B\in K(n-1)$ the balls
$$
\{3L\ :\ L\in K(n,B)\}
$$
are disjoint and contained in $B$ and $3L\subset L^f$.

\item[{\bf(P2)}] For any $n\ge 2$, $B\in K(n-1)$ and any  
$i\in\{1\,\ldots, l_B\}$ the
balls
$$
\{L^f\ :\ L\in K(n,B,i)\}
$$
are disjoint and contained in $B$.

\item[{\bf(P3)}]
For any $n\ge 2$, $B\in K(n-1)$ and  $i\in\{1\,\ldots,l_B\}$ $$
\sum_{L\in K(n,B,i)}V^k(L^f) \ \ge\ c_3\ V^k(B),  $$ where  $c_3
:= \frac{\kappa\,c_1^2}{2\,c_2^2\,10^k} > 0 $ is  an absolute
constant.

\item[{\bf(P4)}] For any $n\ge 2$, $B\in K(n-1)$, any
$i\in\{1\,\ldots,l_B-1\}$ and any $L\in K(n,B,i)$ and $M\in
K(n,B,i+1)$
$$
V^f(M)\le \frac{1}{2}\ V^f(L).
$$

\item[{\bf(P5)}] The number of local sub-levels is defined by
$$
l_B:=\left\{
\begin{array}{lcl}
  \displaystyle\left[\frac{c_2\,\eta}{c_3\,\cH^k(B)}\right]+1 & , &  
\mbox{ if
  $B=B_0 := \K(1)$,}\\[5ex]
  \displaystyle\left[\frac{V^f(B)}{c_3\,V^k(B)}\right]+1 & , & \mbox{ if
  $B\in K(n)$ with $n\ge 2$}
\end{array}
\right.
$$
and satisfies $l_B\ge 2$ for $B\in K(n)$ with $n\ge2$.

\end{enumerate}

\Subsec{The existence of $\K_{\eta}$} \label{kantor1}
In this section we show that it is indeed possible to construct a
Cantor set $\K_{\eta}$ with the desired  properties as discussed
in the previous section.  We will use the notation $$ K_l(n,B) \
:= \ \bigcup_{i=1}^{l}K(n,B,i)  \ . $$ Thus,  $K(n,B)$ is simply
$K_{l_B}(n,B)$.

\demo{Level  $1$} This is defined by
taking the arbitrary ball $B_0$. Thus, $\K(1) := B_0$ and property
(P0) is trivially satisfied.

 \vskip9pt

We proceed  by induction. Assume that the first $(n-1)$ levels
$\K(1), \K(2),\break \ldots , \K(n-1)$ have been constructed. We now
construct the $n$-th level $\K(n)$.

\demo{Level $n$} To construct
this level we construct local levels $K(n,B)$ for each  $B\in
K(n-1)$. Recall, that each local level $K(n,B)$ will consist of
sub-levels $K(n,B,i)$ where $1\leq i \leq l_B$ and $l_B $ is given
by property (P5). Therefore, fix some ball $B\in K(n-1)$ and a
sufficiently small constant $\ve=\ve(B)>0$ which will be
determined later. Let $G$ be sufficiently large so that
\begin{equation}\label{e:019}
\diam(3B_i)<\diam(B_i^f)\qquad\text{ whenever}\qquad i\ge G
\end{equation}

\begin{equation}\label{e:020} \frac{V^k(B_i)}{V^f(B_i)}< \ve \
\frac{V^k(B)}{V^f(B)}\qquad\text{ whenever}\qquad i\ge G
\end{equation}
and
\begin{equation}\label{card}
\left[\frac{V^f(B_i)}{c_3\,V^k(B_i)}\right] \ \geq \ 1
\qquad\text{ whenever}\qquad i\ge G \ ,
\end{equation}
where $c_3$ is the constant appearing in property (P3) above. This
is possible since $x^k/f(x)\to0$ as $x\to0$. Now let ${\cal C}_G
:= \{B_i : i \geq G\}$.   The local level $ K(n,B)$ will be
constructed to be a finite, disjoint sub-collection of ${\cal
C}_G$. Thus, (\ref{e:019})--(\ref{card}) are satisfied for any
ball $B_i$ in $ K(n,B)$. In particular, (\ref{card}) implies that
$l_{B_i} \geq 2 $ and so property (P5) will automatically be
satisfied for balls in $ K(n,B)$.

\demo{Sub-level $1$} With $B$ and $G$ as above, let $\kgb$
denote the collection of balls arising from Lemma~\ref{lem1}.
Note, that in view of (\ref{e:019}) the collection $\kgb$ is a
disjoint collection of balls.  Define the first  sub-level of $
K(n,B)$ to be $\kgb$; that is  $$ K(n,B,1) \  := \ \kgb \ . $$

\noindent By Lemma~\ref{lem1}, it is clear that (P2) and (P3) are
fulfilled for $i=1$. By (\ref{e:019}) and the fact that the balls
in $\kgbf$ are disjoint, we also have that    (P1) is satisfied
within this first sub-level. Clearly, $ K(n,B,1) \subset {\cal
C}_G$.

\demo{Higher sub-levels} To construct higher
sub-levels  we argue by induction. For $l<l_B$,  assume that we
have constructed the  sub-levels $K(n,B,1), \dots\break \dots, K(n,B,l)$
satisfying  properties (P1)--(P4) with $l_B$ replaced by $l$ and
such that $ K_l(n,B) \subset {\cal C}_G$. In view of the latter,
(\ref{e:019})--(\ref{card}) are satisfied for any ball $L$ in $
K_l(n,B)$. In particular, in view of (\ref{card}), for any ball
$L$ in $ K_l(n,B)$ property (P5) is trivially satisfied; i.e. $l_L
\geq 2$.  We now construct the next sub-level $K(n,B,l+1)$.

As every  sub-level of the construction has to be well separated
from the previous ones, we first  verify that there is enough
`space'  left over in $B$  once we have removed the  sub-levels
$K(n,B,1), \dots, K(n,B,l)$ from $B$. More precisely, let
$$A^{(l)} \ := \  \mbox{\small{$\frac{1}{2}$}}B \ \setminus
\bigcup_{L\in K_l(n,B)} \!\!\! 4L \ . $$ We show  that
\begin{equation}\label{e:021}
\cH^k\big(A^{(l)} \big)\ge \frac12\
\cH^k(\mbox{\footnotesize{$\frac{1}{2}$}}B) \ .
\end{equation}

\noindent By construction and the fact that $l <  l_B$,
\begin{eqnarray}&&\label{but}\\  \cH^k(\bigcup_{L\in K_l(n,B)}4L) \!\!\!\!\! & \le
&  \sum_{L\in K_l(n,B)}\cH^k(4L) \nonumber
\\ & \stackrel{(\ref{e:010})}{\,\le\,} &  4^kc_2\sum_{L\in
K_l(n,B)}V^k(L)   \  = \  4^k c_2\sum_{L\in
K_l(n,B)}V^f(L)\,\frac{V^k(L)}{V^f(L)} \nonumber \\  &
\stackrel{(\ref{e:020})}{\,\le\,} & 4^kc_2 \sum_{L\in
K_l(n,B)}V^f(L)\ \ve \ \frac{V^k(B)}{V^f(B)} \nonumber \\ &
\stackrel{(\ref{e:009})}{\,=\,} & 4^kc_2 \ \ve \
\frac{V^k(B)}{V^f(B)}\ \sum_{i=1}^{l}\sum_{L\in K(n,B,i)}V^k(L^f)
\nonumber \\ & \stackrel{(\ref{e:010})}{\,\le\,} & \frac{4^kc_2\
\ve}{c_1}\ \frac{V^k(B)}{V^f(B)}\ \sum_{i=1}^{l}\sum_{L\in
K(n,B,i)}\cH^k(L^f) \nonumber \\ &
  \stackrel{\,\,(P2) \,}{\,\le\,} &
  \frac{4^kc_2\ \ve}{c_1}\
\frac{V^k(B)}{V^f(B)}\ ( l_B -1 ) \ \cH^k(B)  \ .  \nonumber
\end{eqnarray}

\noindent Now, if $B= B_0$  let  $$ \ve = \ve(B_o) \ := \
\frac{1}{2} \left(\frac{c_1}{c_2} \right)^2 \frac{c_3}{2^k \, 4^k}
\ \frac{V^f(B_0)}{\eta} \ .  $$ If $B \neq  B_0$, so that $ B \in
K(n) $ for some $n \geq 2$,   let $ \ve   :=   \ve(B_0) \times
(\eta / V^f(B_0))$ -- a constant independent of $B$, $B_0$ and
$\eta$.
It then follows from (\ref{but}), (P5) and (\ref{e:010})  that $$
\cH^k(\bigcup_{L\in K_l(n,B)}4L) \ \le \ \frac{1}{2} \
\cH^k(\mbox{\footnotesize{$\frac{1}{2}$}}B) \ ,  $$ and this
clearly establishes (\ref{e:021}).

By construction, $K_l(n,B)$ is a finite collection of balls and so
$d_{\min} := \min \{\diam(L): L \in K_l(n,B) \} $ is well defined.
Let $B^{(l)}$ denote a generic ball of diameter $d_{\min}$. At
each point of $A^{(l)} $ place a ball $B^{(l)}$ and denote this
collection by ${\cal A}^{(l)}$. By the $5r$-covering lemma (Lemma
\ref{5r}), there exists a disjoint sub-collection ${\cal G}^{(l)}$
such that $$ A^{(l)} \ \subset \ \bigcup_{B^{(l)} \in {\cal
A}^{(l)} } B^{(l)} \ \subset \ \bigcup_{B^{(l)} \in {\cal G}^{(l)}
} 5 B^{(l)} \ \ \ . $$ The  collection ${\cal G}^{(l)}$ is clearly
contained within $B$  and it  is finite;  the balls are disjoint
and all of the same size. Moreover, by construction
\begin{equation}
  B^{(l)} \cap \bigcup_{L \in K_l(n,B)} 3 L =
\emptyset \hspace{7mm} {\rm for\ any \ }  B^{(l)} \in {\cal
G}^{(l)} \ ; \label{blcoll} \end{equation} i.e.\ the balls in
${\cal G}^{(l)} $ do not intersect any of the $3L$ balls from the
previous sub-levels. It follows that $$ {\cal H}^k (
\bigcup_{B^{(l)} \in {\cal G}^{(l)} } 5 B^{(l)} ) \ \geq \ {\cal
H}^k (A^{(l)} ) \ \stackrel{(\ref{e:021})}{\ \geq \ } \mbox{
\small $\frac12$}  \ \cH^k(\mbox{\footnotesize{$\frac{1}{2}$}}B) \
. $$ On the other hand, since ${\cal G}^{(l)}$ is a disjoint
collection of balls we have that $${\cal H}^k ( \bigcup_{B^{(l)}
\in {\cal G}^{(l)} } 5 B^{(l)} ) \stackrel{(\ref{e:010})}{\ \leq \
} \frac{c_2}{c_1} \ 5^{k} \ \ {\cal H}^k (\bigcup^\circ_{B^{(l)}
\in {\cal G}^{(l)} } B^{(l)} ) \ , $$ and so
\begin{equation}\label{e:025}
{\cal H}^k (\bigcup^\circ_{B^{(l)} \in {\cal G}^{(l)} } B^{(l)}
)\,\ge\,\frac{c_1}{2 c_2 5^k} \ \
\cH^k(\mbox{\footnotesize{$\frac{1}{2}$}}B)\,.
\end{equation}

We are now in the position to construct the $(l+1)$-th sub-level
$K(n,B,\break l+1)$.  To this end, let  $G' \geq G $ be sufficiently
large so that
for every $i\ge G'$
\begin{equation}\label{e:026}
   V^f(B_i) \ \le \  \frac{1}{2}\ \min_{L\in K_l(n,B)}V^f(L)\,.
\end{equation}
We recall that $\{B_i\}$ is the original sequence of balls in
Theorem~\ref{thm3}. The number on the right of (\ref{e:026}) is
well defined and positive as there are only finitely many balls in
$K_l(n,B)$. Furthermore, (\ref{e:026}) is possible since
$\lim_{i\to\infty} \diam(B_i)=0$ and $\lim_{x \to 0}f(x)=0$. Now
to each  ball $B^{(l)} \in {\cal G}^{(l)}$ we apply
Lemma~\ref{lem1} to obtain  a collection $\kgbl$ and  define $$
K(n,B,l+1):=\bigcup_{B^{(l)} \in {\cal G}^{(l)}} \kgbl \ . $$ Note
that since $G' \geq G$, (\ref{e:019})--(\ref{card}) remain valid
and $K(n,B,l+1)\subset {\cal C}_G$. We now verify properties
(P1)--(P5) for this sub-level.

In view of Lemma~\ref{lem1}, for any $B^{(l)} $ in $ {\cal
G}^{(l)}$ the collection $K^f_{G',B^{(l)}}$ is disjoint and
contained within $B^{(l)}$. This together
  with (\ref{e:019}) establishes property (P1)  for balls $L$ in
$\kgbl$. Since the balls $B^{(l)}$ in $ {\cal G}^{(l)}$ are
disjoint and contained within $B$, we have that (P1) is  satisfied
for balls $L$  in $K(n,B,l+1)$. In turn, this together with
(\ref{blcoll}) implies property (P1) for balls  $L$ in
$K_{l+1}(n,B)$. Clearly, the above argument also verifies property
(P2) for balls $L$  in $K(n,B,l+1)$. The following establishes
property (P3) for $i=l+1$: \begin{eqnarray*}
  \sum_{L\in
K(n,B,l+1)} \!\!\! V^k(L^f) \,& = & \, \sum_{B^{(l)} \in {\cal
G}^{(l)} } \ \ \sum_{L \in \kgbl} \!\! V^k(L^f) \\
&\stackrel{(\ref{e:010})}{\,\ge\,} & \frac{1}{c_2}\sum_{B^{(l)}
\in {\cal G}^{(l)} } \ \  \sum_{L \in \kgbl}\cH^k(L^f)\, \\ & & \\
& \stackrel{(\ref{e:014})}{\,\ge\,} & \frac{\kappa}{c_2}
\sum_{B^{(l)} \in {\cal G}^{(l)} } \cH^k(B^{(l)}) \
\stackrel{(\ref{e:025})}{\,\ge\,} \  \frac{\kappa}{c_2}
\frac{c_1}{2 c_2 5^k} \ \
\cH^k(\mbox{\footnotesize{$\frac{1}{2}$}}B)
\\ & &  \\  \ & \stackrel{(\ref{e:010})}{\,\ge\,} & \frac{\kappa  
c_1^2}{2 c_2^2
10^k} \ V^k(B) \ :=  \ c_3 \  V^k(B) \ .
\end{eqnarray*}
  Property
(P4) is trivially satisfied as we have imposed condition
(\ref{e:026}). Finally, in view of (\ref{card}), for any ball $L$
in $K(n,B,l+1)$ property (P5) is satisfied; i.e. $l_{L} \geq 2 $.

  The upshot is that (P1)--(P5) are satisfied up to the
local sub-level\break $K(n,l+1,B)$ and so completes the inductive step.
This establishes the existence of the local level $K(n,B):=
K_{l_B}(n,B)$ for each $B \in K(n-1)$ and thereby the existence of
the $n$-th level $K(n)$.

\Subsec{The measure $\mu$ on $\K_\eta$}
In this section, we  define a probability measure $\mu$ supported
on $\K_\eta$. We will eventually show that the measure satisfies
(\ref{task}). For any ball $L \in K(n)$, we attach a weight
$\mu(L)$ defined recursively as follows.

 For $n =1$, we have that $L= B_0 :=  \K(1)$ and we set  
$\mu(L):=1$.

 For $n \geq 2$, let $L$ be a ball in $K(n)$. By construction,
there is a unique ball $B \in K(n-1)$ such that $L \subset B$. We
set $$ \mu(L)\ := \
\frac{V^f(L)}{\rule{0ex}{2.5ex}\sum\limits_{M\in K(n,B)} \!\!\!
V^f(M)}\ \times \  \mu(B)\ . $$

This procedure thus defines inductively a mass on any
ball appearing  in the construction of $\K_{\eta}$. In fact a lot
more is true; $\mu$ can be further extended to all Borel
subsets $F$ of $\R^k$ to determine $\mu(F)$ so that $\mu$
constructed as above actually defines a measure supported on
$\K_{\eta}$; see Proposition 1.7 \cite{falc}. We state this
formally as a

\demo{Fact}  The probability measure $\mu$ constructed
above is supported on $\K_{\eta}$ and for any Borel subset $F$ of
$\R^k$
$$
\mu(F):= \mu(F \cap  \K_{\eta})  \; = \; \inf\;\sum_{L\in{\cal
C}(F)}\mu(L)  \ ,
$$
where the infimum is taken over all coverings ${\cal C}(F)$ of $F
\cap \K_{\eta}$ by balls  $L \in \bigcup_{n\in\N}K(n)$.

\Subsec{The measure of a  ball in  the Cantor construction}
With $n \geq 2$, the aim of this section is to show that for any
ball $L$ in $ K(n)$  we have that
\begin{equation}\label{e:027}
\mu(L) \ \le \  \frac{V^f(L)}{\eta} \ ;
\end{equation}
i.e.\ (\ref{task}) is satisfied for balls in the Cantor
construction. We start with level $n=2$ and fix a ball $L\in
K(2)=K(2,B_0)$;  recall that $B_0 = \K(1)$. Also, recall that
$B=B^k$ for any ball $B$;  see (\ref{e:008}). By definition, $$
\mu(L) \ := \ \frac{V^f(L)}{\rule{0ex}{2.5ex}\sum\limits_{M\in
K(2,B_0)} \!\!\! V^f(M)}\ \times \ \mu(B_0) \ =\
\frac{V^f(L^k)}{\rule{0ex}{4ex}\sum\limits_{i=1}^{l_{B_0}}\
\sum\limits_{M\in K(2,B_0,i)} \!\!\!  V^f(M^k)} \ . $$ However, $$
\sum\limits_{M\in K(2,B_0,i)} V^f(M^k)\
\stackrel{(\ref{e:009})}{=} \ \sum\limits_{M\in K(2,B_0,i)}
V^k(M^f) \ \stackrel{(P3)}{\geq} \ c_3\ V^k(B_0)\
\stackrel{(\ref{e:010})}{\geq}\  \frac{c_3}{c_2}\ \cH^k(B_0)\,. $$
It now follows from the definition of $l_{B_0}$;  see (P5), that
$$ \mu(L) \ \le \  \frac{c_2\,V^f(L)}{c_3\,\cH^k(B_0)\,l_{B_0}} \
\le \  \frac{V^f(L)}{\eta} \ . $$

\noindent To establish (\ref{e:027}) for general $n$, we proceed
by induction. For $n > 2$, assume that (\ref{e:027}) holds for
balls in $K(n-1)$.  Consider an arbitrary  ball  $L$ in $K(n)$.
Then, $L \in K(n,B)$ for some  $B\in K(n-1)$. By definition and
our induction hypothesis, $$ \mu(L) \ := \
\frac{V^f(L)}{\rule{0ex}{2.5ex}\sum\limits_{M\in K(n,B)} \!\!\!
V^f(M)}\ \times\ \mu(B)\ \le \
\frac{V^f(L)}{\rule{0ex}{2.5ex}\sum\limits_{M\in K(n,B)} \!\!\!
V^f(M)}\times \frac{V^f(B)}{\eta}\ . $$ Thus, (\ref{e:027})
follows on showing  that $$ \sum\limits_{M\in K(n,B)}
V^f(M)\,\,=\,\sum\limits_{M\in K(n,B)} V^f(M^k)\,\ge\, V^f(B) \ .
$$ Well, \begin{eqnarray*}
  \sum\limits_{M\in K(n,B)}
V^f(M^k)\, & = & \,\sum_{i=1}^{l_B}\ \sum\limits_{M\in K(n,B,i)}
V^f(M^k)\  \stackrel{(\ref{e:009})}{=}  \ \sum_{i=1}^{l_B}\
\sum\limits_{M\in K(n,B,i)} V^k(M^f) \\ & & \\
&\stackrel{(P3)}{\ge} & \, c_3 \sum_{i=1}^{l_B}\ V^k(B)\
\stackrel{(P5)}{\ge}\, c_3\ V^k(B)\,\frac{V^f(B)}{c_3\,V^k(B)}\ =\
V^f(B) \ \end{eqnarray*} and so we are done. This completes the
inductive step and thereby establishes (\ref{e:027}) for any $L$
in $K(n)$ with $n \geq 2$.

\Subsec{The measure of an  arbitrary ball}
Set $r_o:=\min\{r(B):B\in K(2)\}$. Take an arbitrary  ball $A$ in
$\R^k$ with $r(A)<r_o$. The aim of this section is to establish
(\ref{task}) for $A$; that is
\begin{equation*}\label{e:028}
\mu(A) \ \ll \ \frac{V^f(A)}{\eta} \ ,
\end{equation*}
were the implied constant is independent of both $A$ and $\eta$.
This will then complete the proof of the Mass Transference
Principle.

We begin by establishing the  following geometric lemma.

\begin{lemma}\label{lem2}
Let $A=B(x_A,r_A)$ and $M=B(x_M,r_M)$ be arbitrary balls such that
$A\cap M\not=\emptyset$ and $A\setminus(cM)\not=\emptyset$ for
some $c\ge3$. Then $r_M\,\le \,r_A$ and $cM\subset 5A$.
\end{lemma}

  {\it Proof.}  Let $z\in A\cap M$. Then $d(x_A,x_M)\le
d(x_A,z)+d(z,x_M)\le r_A+r_M$. Here $d(.,.)$ is the standard
Euclidean metric in $\R^k$. Now take $z\in A\setminus(cM)$. Then
$$ c\,r_M \ \le \  d(x_M,z) \ \le \  d(x_M,x_A)+d(x_A,z) \ < \
r_A+r_M+r_A \ . $$ Hence, $r_M \le \frac{2}{c-1}\ r_A$ and since
$c \geq 3$ we have that $r_M \le r_A$. Now for any $z\in cM$, we
have that \begin{eqnarray*}
  d(x_A,z) \ & \le &  \ d(x_A,x_M)+d(x_M,z) \
\le \ r_A+r_M+ c\,r_M \ =  \  r_A+(1+c)r_M  \\ & & \\  & \le & \
r_A+\frac{2(1+c)}{c-1}r_A \ = \ \Big(3+\frac{4}{c-1}\Big)\ r_A\
\le\ 5\ r_A\ . \end{eqnarray*}
\vglue-20pt\Endproof\vskip12pt

The measure $\mu$ is supported on $\K_\eta$. Thus, without loss of
generality we can assume that $ A \cap \K_\eta \neq \emptyset $;
otherwise $\mu(A) = 0$ and there is nothing to prove.

  We can also assume that for every $n$ large enough $A$ intersects at
  least two balls in  $ K(n)$; since if $B$ is the only ball in  $ K(n)$
which has nonempty intersection with  $A$, then
  $$
  \mu(A) \ \leq \ \mu(B) \ \stackrel{(\ref{e:027})}{\le} \   
\frac{V^f(B)}{\eta} \ \to \ 0
  \hspace{8mm}  {\rm as } \hspace{5mm} n \to \infty  \
  $$
  ($r(B) \to 0$ as $n \to \infty$) and again there is nothing to
  prove. Thus we may assume that there exists a unique  integer $n$
  such that:
\begin{equation}\label{e:029}
\text{$A$ intersects at least 2 balls from $K(n)$}
\end{equation}
and
$$
\text{$A$ intersects only one ball $B$ from $K(n-1)$}.
$$
In view of our  choice of $r_0$ and the fact that $r(A) < r_0$, we
have that $n>2$. Note that since $B$ is the only ball from
$K(n-1)$ which has nonempty intersection with $A$, we trivially
have that $\mu(A)\le\mu(B)$. It follows that we can also assume
that
\begin{equation}
  r(A) \ < \ r(B)  \ .
\label{A<B}
\end{equation}
  Otherwise, since $f$ is increasing
$$
  \mu(A) \ \leq \ \mu(B) \ \stackrel{(\ref{e:027})}{\le} \
   V^f(B)/\eta \ := \ f(\diam(B))/\eta \ \le \  f(\diam(A))/\eta
  \ := \  V^f(A)/\eta  $$
and we are done. Since $K(n,B)$ is a cover for  $A\cap\K_\eta$, we
have that
\begin{equation}\label{e:031}
\mu(A) \ \le \  \sum_{i=1}^{l_B}\ \sum_{L\in K(n,B,i),\ L\cap
A\not=\emptyset} \!\!\!\!\! \mu(L) \ \stackrel{(\ref{e:027})}{\le}
\ \ \sum_{i=1}^{l_B}\ \sum_{L\in K(n,B,i),\ L\cap A\not=\emptyset}
\!\!\!\!\! V^f(L)/\eta\ .
\end{equation}
In order to  estimate the right-hand side  of (\ref{e:031}), we  
consider two
cases:

\medbreak

  \underline{{\em Case}\/} (i):~ Sub-levels $K(n,B,i)$ for
which $$\#\{\,L\in K(n,B,i)\,:\ L\cap A\not=\emptyset\,\} \,=\, 1.$$

\medbreak

  \underline{{\em Case}\/} (ii):~ Sub-levels $K(n,B,i)$
for which $$\#\{\,L\in K(n,B,i)\,:\ L\cap A\not=\emptyset\,\}
\,\ge\, 2.$$
\Enddemo

 Formally,  there is a third case corresponding to those
sub-levels $K(n,B,i)$  for which    $\#\{\,L\in K(n,B,i)\,:\ L\cap
A\not=\emptyset\,\} \,=\, 0$. However, this case is irrelevant
since the contribution to the right-hand side  of (\ref{e:031}) from  
such
sub-levels  is zero.

\demo{Dealing with  Case {\rm (i)}}   Pick a ball $L\in
K(n,B,i)$ such that $L\cap A\not=\emptyset$. By (\ref{e:029}),
there is another ball $M\in K(n,B)$ such that $A\cap
M\not=\emptyset$. By property (P1), $3L$ and $3M$ are disjoint. It
follows that $A\setminus3L\not=\emptyset$. Therefore, by
Lemma~\ref{lem2}, $\diam(L)\le\diam(A)$ and thus
\begin{equation}\label{c}
V^f(L)\ \le \ V^f(A) \ .
\end{equation}
Now, let $K(n,B,i^*)$ denote the first sub-level which has
nonempty intersection with $A$. Thus, $L \cap A = \emptyset$ for
any $L \in K(n,B,i)$ with $i < i^*$ and there exists a unique ball
$L^*$ in  $K(n,B,i^*)$ such that $L^*\cap A \neq \emptyset$. Since
we are in case (i), the internal sum of (\ref{e:031}) consists of
just one summand. It follows, via property (P4) and (\ref{c}),
that
\begin{eqnarray} \sum_{i\,\in\,\text{Case(i)}} \ \sum_{L\in
K(n,B,i),\ L\cap A\not=\emptyset} \!\!\!\!\! V^f(L)/\eta\ & \le &
\ \sum_{i\,\in\,\text{Case(i)}}\ \frac{1}{2^{i-i^*}}\
\frac{V^f(L^*)}{\eta}\,  \label{c1} \\ &\le & \, 2 \;
\frac{V^f(L^*)}{\eta}\,\le\, 2\; \frac{V^f(A)}{\eta}\ .\nonumber
\end{eqnarray}

\demo{Dealing with  Case {\rm (ii)}} Again pick a ball $L\in
K(n,B,i)$ such that $L\cap A\not=\emptyset$. Since we are in case
(ii), there is another ball $M\in K(n,B,i)$ such that $A\cap
M\not=\emptyset$. By property (P2), the balls $L^f$ and $M^f$ are
disjoint. It follows that  $A\setminus L^f\not=\emptyset$. Hence,
by Lemma~\ref{lem2} and property (P1) we have  that
\begin{equation}\label{cc}
L^f\subset 5 A.
\end{equation}
It follows that
\begin{eqnarray} &\!\!\!\!& \label{c2} \\
\sum_{i\,\in\text{ Case(ii)}}  \  \sum_{L\in K(n,B,i),\ L\cap
A\not=\emptyset} \!\!\!\! \frac{V^f(L)}{\eta}  &\!\!
\stackrel{(\ref{e:009})}{=}\!\! &   \sum_{i\,\in\text{ Case(ii)}}\
\sum_{L\in K(n,B,i),\ L\cap A\not=\emptyset} \frac{V^k(L^f)}{\eta}
\nonumber  \\
&  \!\!\!\!& \nonumber \\
&\!\!\stackrel{(\ref{e:010})}{\le}\!\!& \frac{1}{c_1\,
\eta}\sum_{i\,\in\text{ Case (ii)}}\ \sum_{L\in K(n,B,i),\ L\cap
A\not=\emptyset} \!\!\!\!\!\!\! \cH^k(L^f)     \nonumber  \\
& \!\! \!\!& \nonumber \\ & \!\!\stackrel{(P2) \&  (\ref{cc})}{\le}\!\! &
\frac{1}{c_1\,\eta}\sum_{i\,\in\text{ Case(ii)}} \!\!\!\!\!\!
\cH^k(5A) \ \stackrel{(\ref{e:010})}{\le}  \ \frac{5^k c_2
\,V^k(A)\,l_B}{c_1\,\eta}
\nonumber \\
&\!\!\!\!  & \nonumber \\
&\!\!  \stackrel{(P5)}{\le} \!\!& \frac{5^k c_2\,V^k(A)}{c_1\,\eta} \times
\frac{2\,V^f(B)}{c_3\,V^k(B)}\
\nonumber \\
& \!\! \!\!& \nonumber \\
& \!\! \leq\!\! &   \frac{2 \; 5^k c_2}{c_1\,c_3} \times
\frac{V^f(A)}{\eta}\ . \nonumber
\end{eqnarray}
  The last inequality follows
from (\ref{A<B}) and  the fact that the function $x^{-k} f(x)$ is
decreasing.

On combining (\ref{e:031}), (\ref{c1}) and (\ref{c2}) we attain
our goal; i.e.\ $ \mu(A)  \ll  V^f(A)/\eta $.

\vglue-22pt
\phantom{up}

\section{Final comments\label{last}}
\vglue-17pt

\Subsec{A general Mass Transference Principle \label{secmtpg}}
We say that a function $f$ is {\em doubling} if there exists a
constant $\lambda
> 1 $ such that for $x >0$ $$f(2x) \, \leq \, \lambda f(x) \ . $$

Let $(X,d)$ be a locally compact metric space. Let $g$ be a
doubling,  dimension function  and suppose there exist constants
$0<c_1<1<c_2<\infty$ and $r_0 > 0$  such that
\begin{equation*}\label{g}
c_1\ g(\diam(B)) \le \cH^g(B)\le c_2\ g(\diam(B))  \ ,
\end{equation*}
for any ball $B=B(x,r)$ with $x\in X$ and $r\le r_0$. Since $g$ is
doubling, the measure $\cH^g$ is doubling on $X$.  Recall that
$V^g(B):=g(\diam(B))$. Thus, the above condition corresponds to
(\ref{e:010}) in the $\R^k$ setup. Next, given a dimension
function $f$ and a ball $B=B(x,r)$ we define
$$
B^f:=B(x,g^{-1}f(r))\,.
$$
By definition, $B^g(x,r)=B(x,r)$ and
\begin{equation*} \label{e:009n}
V^f(B^g)\,=\,V^g(B^f) \qquad\text{for any ball $B$.}
\end{equation*}
This is an analogue of (\ref{e:009}). In the case $g(x)=x^k$, the
current setup precisely coincides with that of Section \ref{secMTP} in
which $X = \R^k$.    The following result is a natural
generalization of Theorem \ref{thm3} --- the Mass Transference
Principle.

\vglue-18pt
\phantom{up}
\begin{theorem}[A general Mass Transference Principle ]\label{thm3b}
Let $(X,d)$ and $g$ be as above and let $\{B_i\}_{i\in\N}$ be a
sequence of balls in $X$ with $\diam(B_i)\to 0$ as $i\to\infty$.
Let $f$ be a dimension function such that $f(x)/g(x)$ is monotonic
and suppose that for any ball $B$ in $X$
\begin{equation*}\label{e:011b}
\cH^g\big(\/B\cap\limsup_{i\to\infty}B^f_i{}\,\big)=\cH^g(B) \ .
\end{equation*}
Then{\rm ,} for any ball $B$ in $X$
\begin{equation*}\label{e:012b}
\cH^f\big(\/B\cap\limsup_{i\to\infty}B^g_i\,\big)=\cH^f(B) \ .
\end{equation*}
\end{theorem}

The proof of the  general Mass Transference Principle follows on
adapting the  proof of Theorem \ref{thm3}  in the obvious manner.
The property that $\cH^k$ is doubling is used repeatedly in the
proof of Theorem \ref{thm3}. In establishing Theorem~\ref{thm3b},
this property is replaced by  the  assumption  that  $\cH^g$ is
doubling.

  In short, the general Mass Transference Principle allows
us to transfer $\cH^g$-measure theoretic statements for $\limsup$
subsets of $X$ to general $\cH^f$-measure theoretic statements.
Thus, whenever we have a Duffin-Schaeffer type statement with
respect to a  measure $\mu$  comparable to $\cH^g$, we obtain a
general Hausdorff measure theory for free.  For numerous examples
of $\limsup$ sets and associated Khintchine type theorems (the
approximating function $\psi$ is assumed to be monotonic) within
the framework of this section, the reader is referred to
\cite{BDV}.

\Subsec{The Duffin-Schaeffer conjecture revisited} \label{DSR}
Let $\cS_k^* (\psi)$ denote the set of points ${\bf
y}=(y_1,\dots,y_k)\in\I^k$ for which there exist infinitely many
$q\in\N$ and ${\bf p}= (p_1, \dots,p_k) \in \Z^k$ with $(p_1,
\dots,p_k,q) =1$,  such that
\begin{equation*}
  \left|y_i - \frac{p_i}{q}\right|\ <\
\frac{\psi(q)}{q}  \hspace{9mm}   1 \leq i \leq k \ .  \label{1P}
\end{equation*}
  \noindent Here,  we simply ask
that points in $\I^k$ are approximated by distinct rationals
whereas in the definition of $\cS_k(\psi)$ a pairwise co-primeness
condition on the rationals is imposed. For $k=1$, the two sets
coincide. For $k \geq 2$, it is easy to verify that $
m(\cS_k^*(\psi)) = 0 $ if $ \sum   \psi(n)^k   < \infty $. The
complementary divergent result is due to Gallagher \cite{gall}.

\vskip6pt{\scshape Theorem G}  {\em For $k \geq 2${\rm ,} $ \hspace{5mm}
m(\cS_k^*(\psi)) = 1 \ \ \ \ \ {\it if \ } \ \ \ \ \
\sum_{n=1}^{\infty} \ \psi(n)^k \ = \ \infty \ . $}
\vskip6pt

 Notice that the Euler function $\phi$ plays no role in
determining the measure of $\cS_k^*(\psi)$ when $k \geq 2$. This
is unlike  the situation  when considering the measure of the set
$\cS_k(\psi) $;  see Theorem PV (\S\ref{tds}) and Corollary
\ref{cor1}. It is worth mentioning that Gallagher actually obtains
a quantative version of Theorem G.

The Mass Transference Principle together with Theorem G, implies
the following  general statement.

\demo{\scshape Theorem 3}  {\it For} $k \geq 2$, $$  
\cH^f(\cS_k^*(\psi)) = \cH^f(\I^k) \ {\it if} \
 \sum_{n=1}^{\infty} \ f(\psi(n)/n) n^k\break   =   \infty.$$
\Enddemo

  It would be highly desirable to establish a version of
the Mass Transference Principle which allows us to deduce a
quantative Hausdorff measure statement from  a quantative Lebesgue
measure statement.  We hope to investigate this sometime in the
near future.

\demo{Acknowledgments} SV would like to thank Ayesha
and Iona for making him appreciate once again all those
wonderfully simple things around us: ants, pussycats, sticks,
leaves and of course the many imaginary worlds that are often
neglected in adulthood,  especially the world of hobgoblins. VB
would like to thank Tatiana for her help and patience during the
difficult but nevertheless exciting time over the past nine
months.

\end{document}